\documentclass[12pt]{article}
\usepackage{amsmath,subfigure}
\usepackage{amssymb}
\usepackage{epsfig}
\usepackage{eepic}
\usepackage{graphicx}
\usepackage{float}
\usepackage{latexsym,graphics,color}
\usepackage{dsfont}
\newtheorem{theorem}{Theorem}[section]
\newtheorem{proposition}[theorem]{Proposition}
\newtheorem{lemma}[theorem]{Lemma}

\newtheorem{definition}[theorem]{Definition}

\newtheorem{assumption}[theorem]{Assumption}

\newtheorem{remark}{Remark}[section]

\def\sqr#1#2{{\vcenter{\vbox{\hrule height.#2pt\hbox{\vrule width.#2pt height#1pt \kern#1pt\vrule width.#2pt}\hrule height.#2pt}}}}

\newcommand{\R}{{\mathsf{I\!R}}}
\allowdisplaybreaks

\begin{document}
                   \begin{center}
   {\large \bf Large-time behavior of a two-scale semilinear reaction-diffusion system for concrete sulfatation}\\[1cm]
                   \end{center}

                   \begin{center}
                  {\sc Toyohiko Aiki$^{\rm \ddag}$ 
          and Adrian Muntean$^{\rm \star}$  }\\
                    \vspace{0.3cm}

 $^{\rm \ddag}$Department of Mathematics,
Faculty of Science,
Japan Women's University \\
2-8-1 Mejirodai, Bunkyo-ku, Tokyo
112-8681, Japan. e-mail: aikit@fc.jwu.ac.jp

$^{\rm \star}$CASA - Centre for Analysis,
Scientific computing and Applications, Institute for Complex Molecular Systems (ICMS), Department of Mathematics and
Computer Science, Eindhoven University of Technology, The Netherlands.  e-mail: a.muntean@tue.nl

        \end{center}
       \vspace{1cm}
\begin{abstract} We study the large-time behavior of (weak) solutions to a two-scale reaction-diffusion system coupled with a nonlinear ordinary differential equations modeling the partly dissipative corrosion of concrete (/cement)-based materials with sulfates. We prove that as $t\to\infty$ the solution to the original two-scale system converges to the corresponding two-scale stationary system. To obtain the main result we make use essentially of the theory of evolution equations governed by subdifferential operators of time-dependent convex functions developed combined with a series of two-scale energy-like time-independent estimates.
\end{abstract}

 \section{Introduction}
 \subsection{Organization of the paper}
 
 We study the large time behavior of a two-scale reaction-diffusion system modeling the evolution of the sulfatation reaction in concrete-based materials; see \cite{Tasnim2scale} for a rigorous derivation of the system by periodic homogenization [a direct application of two-scale convergence principles  \cite{Allaire,Nguestseng} and multiscale analysis of PDEs posed in perforated domains \cite{Hornung_Jager}]. To fix ideas, let us only mention here that the sulfatation reaction attacks aggressively unsaturated porous media, where the ${\rm H_2S}$ air-water transfer and bacteria interplay together in the presence of heat. This is precisely the case of most sewer pipes or of marble monuments in countries like Brazil, Japan, USA, Italy,  etc.; see e.g. \cite{Rosen,Natalini} and references cited therein. More engineering details on this scenario can be found e.g. in \cite{Gleize}.

To show  that as $t\to\infty$ the solution to the original two-scale reaction-diffusion system converges to the corresponding two-scale stationary system, we proceed as follows:

 In the subsequent sections, we present the setting of the two-scale model equations (Section \ref{Model}), give a brief outlook on the literature that inspired us to working in such framework of multiple spatial scales (Section \ref{lit}), and finally, we delimitate the aim of the paper (Section \ref{aim}). The main technical assumptions behind our results together with the weak solvability of the problem are collected in Section \ref{setting} and Section \ref{wp}, respectively. The bulk of the paper is  Section \ref{LARGE} -- the proof of our main main result - Theorem \ref{large} - a characterization of the solution behavior at large times.

\subsection{Two-scale model equations}\label{Model}
Let us  consider $\Omega$ and $Y$ to be connected and bounded domains in 
$\R^3$, 
$\Gamma_i \subset \partial Y$, $i = 1, 2, 3$, $\partial Y = \Gamma_1 \cup \Gamma_2 \cup \Gamma_3$, and 
$\Gamma_1$, $\Gamma_2$ and $\Gamma_3$  are disjoint. 
Also, $\partial \Omega = \Gamma_D \cup \Gamma_N$ and $\Gamma_D \cap \Gamma_N =\emptyset$. 

Simplifying the scenario analyzed in \cite{two-scale1},  we consider here the following system of partial differential equations coupled with one ordinary differential equation for $T > 0$:
\begin{eqnarray}\label{main}
&& \partial_t w_1 - \nabla_y\cdot (d_1 \nabla_y w_1)= - \psi(w_1 - \gamma w_2) \quad  
\text{ in }
(0,T)\times  \Omega \times Y, \label{pd1}\\
&&  \partial_t  w_2 - \nabla_y \cdot (d_2 \nabla_y w_2)=  \psi(w_1 -  \gamma w_2)
 \quad \text{ in } (0,T)\times\Omega \times Y, \label{pd2}\\
&& \partial_t  w_3 - \nabla\cdot (d_3 \nabla w_3)= -\alpha\int\limits_{\Gamma_2}\big( h w_3 - w_2\big)d\gamma_y
   \quad \text{ in } (0,T)\times\Omega, \label{pd3}\\
&&  \partial_t w_4 = \eta (w_1, w_4)  \quad
  \text{ on } (0,T)\times\Omega\times\Gamma_1.\label{od1}
\end{eqnarray}
 The system is equipped with the initial conditions
\begin{equation}
 w_i(0, \cdot, \cdot) = w_{i0}  \text{ in } \Omega\times Y,  i = 1,2,  
w_3(0, \cdot) = w_{30}   \text{ in } \Omega,  w_4(0, \cdot, \cdot) = w_{40}  \text{ on } \Omega \times\Gamma_1, \label{ic}
\end{equation}
while the boundary conditions are
\begin{equation}
\begin{cases}
d_1 \nabla_y w_1 \cdot \nu_y =  -\eta (w_1, w_4)   \text{ on } (0,T)\times\Omega \times\Gamma_1, \\
d_1 \nabla_y w_1 \cdot \nu_y =  0  \text{ on } (0,T)\times \Omega 
  \times (\Gamma_2 \cup \Gamma_3), \\
d_2 \nabla_y w_2 \cdot \nu_y =  0  \text{ on } (0,T) 
\times \Omega \times (\Gamma_1 \cup \Gamma_3), \\
d_2 \nabla_y w_2 \cdot \nu_y =\alpha(h w_3-w_2\big)  
 \text{ on }  (0,T)\times\Omega \times\Gamma_2,\\
d_3 \nabla w_3 \cdot \nu(x) =  0    \text{ on } (0,T)\times\Gamma_N, \quad 
\;w_3 = w_3^D  
\text{ on } (0,T)\times\Gamma_D, \label{main_bc}
\end{cases}
\end{equation}
where  $w_1 =w_1(t,x,y)$ denotes the concentration of H$_2$SO$_4$ in $(0,T)\times\Omega\times Y$, $w_2 = w_2(t,x,y)$ the concentration of H$_2$S aqueous species in $(0,T)\times\Omega\times Y$,
$w_3 = w_3(t,x)$ the concentration of H$_2$S gaseous species in $(0,T)\times\Omega$ and $w_4 =w_4(t,x,y)$ of $gypsum$ concentration on $(0,T)\times\Omega\times\Gamma_1$ and $\eta$ is the reaction rate of $gypsum$. 
$\nabla$ without subscript denotes the differentiation w.r.t. macroscopic variable $x$, while $\nabla_y$ is the respective differential operators w.r.t. 
the micro-variable $y$, $\nu$ and $\nu_y$ are outward normal vectors on 
$\partial \Omega$ and $\partial Y$, respectively. 

By  $\alpha$ we denote the rate of the reaction taking place on the interface $\Gamma_2$, $h$ is Henry's constant (see \cite{Henry} for an extensive review on Henry's law), $d_1$, $d_2$ and $d_3$ are diffusion coefficients,  $\psi$ is a continuous function on $\R$ and $\gamma$ is a positive constant. The microscale and macroscale are 
 connected together via the right-hand side of $\eqref{pd3}$ and via the
{\em micro-macro} boundary condition $\eqref{main_bc}_4$. 

\subsection{Comments on related PDE systems with multiple scales structure}\label{lit}

 Promoted initatially by  G. I. Barenblatt and his co-workers  (see e.g. \cite{Barenblatt} and references cited therein), the  research on fissured media-like equations in particular and on double-porosity models in general has attracted an increasing interest on the mathematical analysis side  especially on what pseudo-parabolic equations and nonlinear parabolic equations posed on multiple spatial scales are concerned; the focus being mostly on well-posedness aspects; see the original work of  R.E Showalter and collaborators \cite{Show_MB,Show_Hornung,Show_walk, Show_Cook} which was a source of inspiration for later developments by M. B\"ohm and S. A. Meier (\cite{Meier_Bohm,sebam_PhD,CRAS,fast}) and J. Escher and D. Treutler  (\cite{Daniela1,Daniela2}).

 Mathematical tools employed range from a fine use of  strong solutions (exploiting the semigroup structure of the problem), fixed-point arguments in Bochner spaces, energy methods for parabolic equations as well as weak-convergence type methods (particularly the two-scale convergence put in the periodic homogenization context). 

For a classification of two-scale PDE systems based on the used  {\em micro-macro transmission condition}, we refer the reader to the chapter written by R. E. Showalter in U. Hornung's book \cite{hornung}. Furthermore, the reader will discover therein that the concept of balance equations on two-scales (or on a {\em distributed array of microstructures}) can be used as a stand alone modeling tool, not necessarily in the context of  averaging techniques.

 \subsection{Aim of this paper}\label{aim}
 
Very much in line with older results by Friedman, Knabner,  and  Tzavaras (compare \cite{Friedman,Knabner}),  the interest of this paper lies on the large-time asymptotics of the two-scale system (\ref{pd1})--(\ref{od1}) endowed with the initial conditions (\ref{ic}) and the  boundary conditions  (\ref{main_bc}). 

In this special context of multiple spatial scales, we need to cope with two main specific difficulties:
\begin{itemize}
\item[(i)] Due to the presence of the microscale (here denoted by $Y$)  and evolution equations posed at that level, memory effects are inherently present; see e.g. \cite{Nenad}. The question is here twofold: How strong are such memory effects and to which extent can they affect the lifespan of the coupled PDE system? 
 
\item[(ii)] As micro-macro transmission condition we impose a nonlinear Henry's law, therefore particular care is needed while treating two-scale traces of Sobolev functions; see e.g. \cite{MunteanNeuss2010}.
\end{itemize}  
It is worth noting that the large-time behavior is not only the most interesting mathematical question that one would pose at this stage, but also it  is the most relevant one from the practical point of view -- an estimate on the lifespan of the material [forced to confront, for instance,  evolving free boundaries \cite{Aiki_CPAA, Aiki_IFB}, potential  clogging of the pores, and self-healing \cite{Vermolen}] is the holy grail of the materials science.  

For our problem (\ref{pd1})--(\ref{main_bc}), we basically show that the large-time behavior of the active concentrations is well described by the solution of the corresponding stationary system; see Theorem \ref{large}. In the rest of the paper, we prepare a suitable mathematical framework and  then prove the large-time asymptotics for this  multiscale reaction-diffusion scenario in  a rigorous manner.   

\section{Assumptions and  main results}\label{setting}

To keep notation simple, we put 
$$ X := \{z \in H^1(\Omega): z=0 \mbox{ on } \Gamma_D \}, 
 H := L^2(\Omega \times Y) \times L^2(\Omega \times Y)  \times L^2(\Omega), $$
$$ V := L^2(\Omega; H^1(Y)) \times  L^2(\Omega; H^1(Y)) \times X, 
K(T) := W^{1,2}(0,T; L^2(\Omega \times \Gamma_1)) \mbox{ for } T > 0, $$
$$ (u,v)_H := (u_1,v_1)_{L^2(\Omega \times Y)} + 
 \gamma (u_2,v_2)_{L^2(\Omega \times Y)} + \gamma h(u_3,v_3)_{L^2(\Omega)} $$
for $u:= (u_1, u_2, u_3)$, $v := (v_1, v_2, v_3) \in H$.

\begin{assumption}\label{assump}
(A1) For i = 1, 2, 
$d_i\in L^\infty(\Omega\times Y)$ and $d_3\in L^\infty(\Omega)$ such that 
$d_i(x,y) \geq d^0_{i}$  for a.e. $(x,y)\in \Omega\times Y$, where 
$d_i^0$ is a positive constant, and 
$d_3(x) \geq d^0_{3}$ for 
a.e. $x \in \Omega$, where $d_3^0$ is a positive constant.

(A2) $\eta(r_1, r_2):= R(r_1)Q(r_2)$  for $r_1, r_2 \in \R$,  
where $R$ and $Q$ are locally Lipschitz continuous functions such that  
$R'\geq 0$ and $Q'\leq 0$ a.e. on $\R$ and 
   \begin{eqnarray}
R(r_1):=
\left \{ \begin{array}{cl}
\mbox{positive} & \mbox{ if } r_1 > 0, \\ 
0 & \mbox{ otherwise, }
 \end{array} \right.
\quad\quad
Q(r_2):=
\left \{ \begin{array}{cc}
\mbox{positive} & \mbox{ if } r_2 < \beta_{max},\\ 
0  & \mbox{otherwise, } 
 \end{array} \right. 
\nonumber 
\end{eqnarray}
where $\beta_{max}$ is a positive constant. 
 
(A3)  The function $\psi$ is increasing and locally Lipschitz continuous 
on $\R$ with $\psi(0)= 0$.

(A4) $w^D_{3}\in L^2_{loc}([0,\infty);H^2(\Omega))\cap 
W^{1,2}_{loc}([0, \infty); L^2(\Omega))\cap L_+^\infty((0,\infty)\times\Omega)$ with $\nabla w_3^D \cdot \nu = 0$ on $(0,\infty) \times \Gamma_N$.

(A5)  $w_{i0}\in L^2(\Omega; H^1(Y))\cap L^\infty_+(\Omega\times Y)$ 
for $i = 1, 2$, 
$w_{30}\in H^1(\Omega)\cap L^\infty_+(\Omega)$, 
 $w_{30} - w^D_3(0,\cdot) \in X$, and  $w_{40}\in L^\infty_+(\Omega\times \Gamma_1)$.
\end{assumption}
Note that in (A4) and (A5) we define $L_+^{\infty}(\Omega'): = L^{\infty}(\Omega') \cap \{u:  u \geq 0 \mbox{ on } \Omega'\}$ for a domain $\Omega'$.

Next, we denote the two-scale problem \eqref{pd1}--\eqref{main_bc} by TP$(R, Q, \psi)$ and give a definition of a solution to TP$(R, Q, \psi)$ as follows:
\begin{definition}\label{def}
We call the multiplet $(w_1, w_2, w_3, w_4)$ a solution to TP$(R, Q, \psi)$  on $[0,T]$, $T > 0$, if (S1) $\sim$ (S5) hold: 
\\ 
(S1) $w_1, w_2 \in W^{1,2}(0,T; L^2(\Omega \times Y)) \cap L^{\infty}(0,T; L^2(\Omega; H^1(Y))) 
\cap L^\infty((0,T)\times\Omega\times Y)$, 
$w_3  \in W^{1,2}(0,T; L^2(\Omega)) \cap L^{\infty}((0,T) \times \Omega)$, 
$w_3 - w_3^D \in  L^\infty(0,T; X)$, 
 $w_4 \in W^{1,2}(0,T; L^2(\Omega\times\Gamma_1)) \cap L^\infty((0,T)\times\Omega\times \Gamma_1)$. 
$$ \int_{\Omega \times Y} \partial_t w_1 v_1  dxdy 
+ \int_{\Omega \times Y} d_1 \nabla_y w_1 \cdot \nabla_y  v_1 dxdy 
 + \int_{\Omega \times \Gamma_1} Q(w_4) R(w_1) v_1 dxd\gamma_y 
\leqno{(S2)} $$
$$=   - \int_{\Omega \times Y} \psi (w_1 - \gamma w_2)    v_1 dxdy 
\quad \mbox{ for } v_1 \in L^2(\Omega; H^1(Y)) \mbox{ a.e. on } [0,T]. 
$$
$$
 \int_{\Omega \times Y} \partial_t w_2  v_2 dxdy 
+ \int_{\Omega \times Y} d_2 \nabla_y w_2 \cdot \nabla_y v_2 dxdy 
 - \alpha  \int_{\Omega \times \Gamma_2} (h w_3 - w_2) v_2 dxd\gamma_y 
\leqno{(S3)} $$
$$ =  
 \int_{\Omega \times Y} \psi(w_1 -  \gamma w_2)  v_2 dxdy 
 \quad 
\mbox{ for } v_2 \in L^2(\Omega; H^1(Y))  \mbox{ a.e. on } [0,T]. $$
$$ \int_{\Omega} \partial_t w_3  v_3 dx 
+ \int_{\Omega } d_3 \nabla w_3 \cdot \nabla v_3 dx 
=  - \alpha  \int_{\Omega \times \Gamma_2} (h w_3  - w_2) v_3 
 dxd\gamma_y  
  \mbox{ for } v_3 \in X \mbox{ a.e. on } [0,T]. \leqno{(S4)} 
$$
(S5) (\ref{od1})  and  (\ref{ic}) hold. 

Moreover, the multiplet $(w_1, w_2, w_3, w_4)$ is called a solution of TP$(R, Q, \psi)$ 
on $[0,\infty)$, if it is a solution of TP$(R, Q, \psi)$ on $[0,T]$ for any $T > 0$. 
\end{definition}

These two Theorems are concerned with the well-posedness of TP$(R, Q, \psi)$. 
\begin{theorem}\label{Uniqueness}(Uniqueness)
Assume (A1) $\sim$ (A5), then there exists at most one solution of TP$(R, Q, \psi)$. 
\end{theorem}

\begin{theorem} \label{Exist}
(Global existence of solutions to TP$(R, Q, \psi)$) 
Assume (A1) $\sim$ (A5), then there exists a solution $(w_1,w_2, w_3, w_4)$ of 
TP$(R, Q, \psi)$ on $[0,\infty)$. Moreover, it holds that 

(i) $w_1(t),w_2(t)\geq0$ a.e. in ${\Omega\times Y}$, $w_3(t) \geq 0$ a.e. in ${\Omega}$ and
 $w_4(t)\geq0$ a.e. on $\Omega\times\Gamma_1$ for $t \geq 0$. 

(ii) $w_1(t)\leq M_1$, $w_2(t)\leq M_2$ a.e. in ${\Omega\times Y}$,
 $w_3(t)\leq M_3$ a.e. in ${\Omega}$ and
$w_4(t)\leq M_4$ a.e. on $\Omega\times\Gamma_1$ for $t \geq 0$, 
 where 
$M_1 =$ \\ $\max\{ |w_{10}|_{L^{\infty}(\Omega \times Y)},  
 \gamma   |w_{20}|_{L^{\infty}(\Omega \times Y)}, 
\gamma h  |w_{30}|_{L^{\infty}(\Omega)},  
   \gamma h |w_{3}^D|_{L^{\infty}((0,\infty) \times \Omega)}\}$, 
$M_2 = \frac{1}{\gamma} M_1$, $M_3 = \frac{1}{h} M_2$, 
and $M_4 = \max\{ \beta_{max}, |w_{40}|_{L^{\infty}(\Omega \times \Gamma_1)}\}$.\end{theorem}

\vskip 12pt
\begin{remark}
In \cite{two-scale1},  we consider the system consisting in  
the  equations 
 (\ref{pd3}) and (\ref{od1}) and prove its  well-posedness (see Theorem 4.3 and Theorem 4.4 in \cite{two-scale1}): 
\begin{eqnarray*}
&& \partial_t w_1 - \nabla_y\cdot (d_1 \nabla_y w_1)= - f_1(w_1) + f_2(w_2) 
\quad  
\text{ in } (0,T)\times  \Omega \times Y, \\
&&  \partial_t  w_2 - \nabla_y \cdot (d_2 \nabla_y w_2)=  f_1(w_1) -  f_2(w_2)
 \quad \text{ in } (0,T)\times\Omega \times Y, 
\end{eqnarray*}
where $f_1$ and $f_2$ are continuous functions on $\R$. 
It is easy to see that Theorems \ref{Uniqueness} and \ref{Exist} can cover the well-posedness of TP$(R, Q, f_1-f_2)$ in case $f_1(r) = a_1[r]^+$ and $f_2(r) = a_2[r]^+$ for 
$r \in \R$, where $a_1$ and $a_2$ are positive constants. In fact, 
Let $\psi(r) = a_1(r)$ for $r \in \R$ and  $\gamma = \frac{a_2}{a_1}$.  
Then, by Theorem \ref{Exist},   the solution $(w_1, w_2, w_3, w_4)$ 
of TP$(R, Q, \psi)$ satisfies that
\begin{itemize}
 \item  $w_1$ and $w_2$ are nonnegative;
\item $\psi(w_1 - \gamma w_2) = f_1(w_1) - f_2(w_2)$. 
\end{itemize}
\end{remark}

To be able to study  the large time behavior of the solution,  we need the additional  condition 
(A6) on the boundary data.

(A6)  $w_3^D \in L^{\infty}(0,\infty; H^1(\Omega))$, 
$\partial_t  w_{3}^D  \in L^1((0,\infty) \times \Omega) \cap 
L^2((0,\infty) \times \Omega) $, 
$\partial_t w_{3}^D - \nabla d_3 \nabla w_{3}^D \in L^{\infty}(0,\infty; L^1(\Omega))$, 
$\partial_t (\partial_t w_{3}^D - \nabla d_3 \nabla w_{3}^D) \in 
L^1((0,\infty) \times \Omega)$.  

\vskip 12pt
Clearly, under (A6) we have: 
$w_3^D(t) \to w_{3\infty}^D$ in $L^2(\Omega)$ and weakly in $H^1(\Omega)$ 
as $t \to \infty$. 
Also, in order to give a statement on the large time behavior 
we introduce the following stationary problem 
SP($w_{4\infty}$, $w_{3\infty}^D$) for given functions 
$w_{4\infty}$ and $w_{3\infty}^D$. In this problem unknown functions are $w_{1\infty}$, $w_{2\infty}$ and $w_{3\infty}$ such that 
\begin{eqnarray*} 
 &&  - \nabla_y \cdot (d_1 \nabla_y w_{1\infty}) 
= - \psi(w_{1\infty} - \gamma w_{2\infty}) \quad 
            \text{ in } \Omega \times Y,  \\
&&  - \nabla_y \cdot (d_2 \nabla_y w_{2\infty}) = 
 \psi(w_{1\infty}  -  \gamma w_{2\infty}) \quad 
       \text{ in }  \Omega \times Y, \\
&& - \nabla\cdot (d_3 \nabla w_{3\infty}) = - \alpha \int_{\Gamma_2}
     \big(h w_{3\infty} - w_{2\infty} \big)d\gamma_y 
   \quad \text{ in } \Omega, \\
&& d_1 \nabla_y w_{1\infty} \cdot \nu_y =  -\eta (w_1, w_{4\infty}) 
  \text{ on }  \Omega \times\Gamma_1, \quad 
 d_1 \nabla_y w_{1\infty} \cdot \nu_y =  0  
         \text{ on } \Omega \times (\Gamma_2 \cup \Gamma_3), \\   
&& d_2 \nabla_y w_{2\infty} \cdot \nu_y =  0  \text{ on } 
       \Omega \times (\Gamma_1 \cup \Gamma_3), \quad 
 d_2 \nabla_y w_{2\infty} \cdot \nu_y = \alpha(h w_{3\infty} - 
   w_{2\infty} \big)  \text{ on } \Omega \times\Gamma_2,\\
&& d_3 \nabla w_{3\infty} \cdot \nu(x) =  0    \text{ on } \Gamma_N, \quad
       w_{3\infty} = w_{3\infty}^D  \text{ on } \Gamma_D.  
\end{eqnarray*}

\begin{definition}\label{def_SP}
We say that the triplet $(w_{1\infty}, w_{2\infty}, w_{3\infty})$ is a solution of SP($w_{4\infty}$, $w_{3\infty}^D$), if the following conditions hold: 
$w_{i\infty} \in L^2(\Omega; H^1(Y)) \cap L^{\infty}(\Omega \times Y)$, $i = 1,2$, $w_{3\infty} - w_{3\infty}^D \in X \cap L^{\infty}(\Omega)$, 
\begin{eqnarray*}
&  & \int_{\Omega \times Y} d_1 \nabla_y w_{1\infty} \cdot \nabla_y v_1 dxdy 
    + \int_{\Omega \times \Gamma_1} Q(w_{4\infty}) R(w_{1\infty}) 
 v_1 dxd\gamma_y
\\
& = & - \int_{\Omega \times Y} \psi(w_{1\infty} -\gamma w_{2\infty}) v_1 dxdy 
      \quad \mbox{ for } v_1 \in L^2(\Omega, H^1(Y)), 
\end{eqnarray*}
\begin{eqnarray*}
&  & \int_{\Omega \times Y} d_2 \nabla_y w_{2\infty} \cdot \nabla_y v_2 dxdy 
  - \alpha \int_{\Omega \times \Gamma_2}  (h w_{3\infty}  - 
    w_{2\infty}) v_2  dxd\gamma_y
\\
& = &  \int_{\Omega \times Y} \psi(w_{1\infty} -\gamma w_{2\infty}) v_2 dxdy 
      \quad \mbox{ for } v_2 \in L^2(\Omega, H^1(Y)), 
\end{eqnarray*}
$$ \int_{\Omega} d_3 \nabla w_{3\infty} \cdot \nabla v_3 dx 
  + \alpha \int_{\Omega \times \Gamma_2}  (hw_{3\infty} - w_{2\infty}) v_3 
  dxd\gamma_y 
 =  0 \quad  \mbox{ for } v_3 \in X. $$
 \end{definition}

\begin{theorem}\label{large}(Large-time behavior)
Assume (A1) $\sim$ (A6)  and 
let $(w_1, w_2, w_3, w_4)$ be a solution  of TP$(R, Q, \psi)$ on $[0,\infty)$. 
Then it holds that 
$w_4(t) \to w_{4 \infty}$ in $L^1(\Omega \times \Gamma_1)$ as $t \to \infty$ 
for some $w_{4\infty} \in L^{\infty}(\Omega \times \Gamma_1)$, 
and $\partial_t w_4 \in L^{1}((0,\infty) \times \Omega \times \Gamma_1)$ and 
there exists a sequence $\{t_n\}$ with $t_n \to \infty$ as $n \to \infty$ such that 
$$ w(t_n) \to w_{\infty} \mbox{ weakly in } H \mbox{ as } n \to \infty, $$
for some $w_{\infty} \in H$, and $w_{\infty}$ is a solution of  
SP($w_{4\infty}$, $w_{3\infty}^D$), where  $w = (w_1, w_2, w_3)$. 

Moreover, if $(\psi(r) - \psi(r'))(r - r') \geq \mu |r - r'|^{p+1}$ for $r, r' \in \R$, where $\mu > 0$ and $p \geq 1$, then 
SP($w_{4\infty}$, $w_{3\infty}^D$) has at most one solution 
and   
\begin{equation} 
 w(t) \to w_{\infty} \mbox{ in } H \mbox{ as } t \to \infty. \label{as3}
\end{equation} 
\end{theorem}

\section{Well-posedness of TP$(R, Q, \psi)$}\label{wp}
The aim of this section is to show the existence and uniqueness of a solution to 
TP$(R, Q, \psi)$ on $[0,T]$ for any $T > 0$. 

First, we consider the following auxiliary problem AP($w_4$) for given 
$w_4 \in  K(T)$: 
\begin{eqnarray} 
&& \partial_t w_1 - \nabla_y\cdot (d_1 \nabla_y w_1)= - \psi(w_1 - \gamma w_2) \quad  
\text{ in }
(0,T)\times  \Omega \times Y, \label{pd1b}\\
&&  \partial_t  w_2 - \nabla_y \cdot (d_2 \nabla_y w_2)=  \psi(w_1 -  \gamma w_2)
 \quad \text{ in } (0,T)\times\Omega \times Y, \label{pd2b}\\
&& \partial_t  w_3 - \nabla\cdot (d_3 \nabla w_3)= -\alpha\int\limits_{\Gamma_2}\big(h w_3-w_2\big)d\gamma_y
   \quad \text{ in } (0,T)\times\Omega, \label{pd3b} \\
&&  w_i(0) = w_{i0}    \text{ in } 
    \Omega\times Y \mbox{ for } i = 1,2, \mbox{ and } w_3(0) = w_3^0
   \text{ in } \Omega, \\
&& d_1 \nabla_y w_1 \cdot \nu_y =  -\eta (w_1, w_4) \quad  \text{ on } (0,T)\times\Omega \times\Gamma_1, \\
&& d_1 \nabla_y w_1 \cdot \nu_y =  0  \quad \text{ on } 
    (0,T)\times\Omega \times ( \Gamma_2 \cup \Gamma_3), \\
&& d_2 \nabla_y w_2 \cdot \nu_y =  0 \quad \text{ on } 
  (0,T)\times\Omega \times ( \Gamma_1 \cup \Gamma_3),  \\
&& d_2 \nabla_y w_2 \cdot \nu_y =\alpha(h w_3-w_2\big)  \quad 
 \text{ on }  (0,T)\times\Omega \times\Gamma_2,\\
&& d_3 \nabla w_3 \cdot \nu(x) =  0  \text{ on } (0,T)\times\Gamma_N, w_3 = w_3^D  \text{ on } (0,T)\times\Gamma_D.  \label{bdb}
\end{eqnarray}

From now on, we solve the above problem AP($w_4$) by using the theory of evolution equations governed by subdifferential operators of time-dependent convex functions (see \cite{Yamada} and\cite{Ke(Chiba)}). To apply this theory,  we first define 
a function $\varphi^t(w_4, \cdot)$ on $H$ for $t \in [0,T]$, 
$w_{3}^D$ satisfying (A4) and 
$w_4 \in  L^2(\Omega \times \Gamma_1)$ as follows: 
\begin{eqnarray*}
& & \varphi^t(w_4, u) \\
& = & 
\left\{ \begin{array}{l} \displaystyle
\frac{1}{2} \int_{\Omega \times Y} d_1 |\nabla_y u_1|^2 dxdy 
+  \int_{\Omega \times \Gamma_1} Q(w_4)  \hat{R}(u_1) dxd\gamma_y 
+ \frac{\gamma}{2} \int_{\Omega \times Y} d_2 |\nabla_y u_2|^2 dxdy \\
\displaystyle
+  \int_{\Omega \times Y} \hat{\psi} (u_1 - \gamma u_2) dxdy  
 + \frac{\gamma \alpha}{2} \int_{\Omega \times \Gamma_2} 
         |h(u_3 + w_3^D(t)) - u_2|^2 dxd\gamma_y 
 \\ 
\displaystyle
  + \frac{\gamma h}{2} \int_{\Omega} d_3|\nabla u_3|^2 dx \quad 
\mbox{ if } u = (u_1, u_2, u_3) \in V, \\
\infty \quad \quad  \mbox{ otherwise, }
\end{array} \right. 
\end{eqnarray*}
where $\hat{R}$ and $\hat{\psi}$ are primitives of $R$ and $\psi$ with 
$\hat{R}(0) = 0$ and $\hat{\psi}(0) = 0$, respectively.
Moreover, we can prove the following Lemma in a straightforward manner:

\begin{lemma} \label{lem1}
Let $t \in [0,T]$. 
If (A1) $\sim$ (A4) hold,  
$Q$ is Lipschitz continuous and bounded on $\R$, 
$R$ and $\psi$ are Lipschitz continuous on $\R$  and 
 $w_4 \in K(T)$, then $\varphi^t(w_4, \cdot)$ is proper, l.s.c. and convex on $H$ for $t \in [0,T]$ and $D(\varphi^t(w_4, \cdot)) = V$ and 
$\partial \varphi^t(w_4, u)$ is single valued, 
where  $D(\varphi^t(w_4, \cdot))$ denotes the effective domain of 
$\varphi^t(w_4, \cdot)$. Moreover, 
$u \in H$ and $u^* = \partial \varphi^t(w_4,u)$ is equivalent to 
$u^* = (u_1^*, u_2^*, u_3^*) \in H$ and 
\begin{eqnarray*}
(u_1^*, v_1)_{L^2(\Omega \times Y)} 
& = & \int_{\Omega \times Y} d_1 \nabla_y u_1 \cdot \nabla_y v_1 dxdy 
  + \int_{\Omega \times \Gamma_1} Q(w_4) R(u_1) v_1 dxd\gamma_y
\\
& & + \int_{\Omega \times Y} \psi(u_1 -\gamma u_2) v_1 dxdy 
      \quad \mbox{ for } v_1 \in L^2(\Omega, H^1(Y)), \\
(u_2^*, v_2)_{L^2(\Omega \times Y)} 
& = & \int_{\Omega \times Y} d_2 \nabla_y u_2 \cdot \nabla_y v_2 dxdy 
  - \alpha \int_{\Omega \times \Gamma_2}  (h(u_3 + w_3^D(t)) - u_2) v_2 
  dxd\gamma_y
\\
& & - \int_{\Omega \times Y} \psi(u_1 -\gamma u_2) v_2 dxdy 
      \quad \mbox{ for } v_2 \in L^2(\Omega, H^1(Y)), 
\end{eqnarray*}   
$$
(u_3^*, v_3)_{L^2(\Omega)} 
 =  \int_{\Omega} d_3 \nabla u_3 \cdot \nabla v_3 dx 
  + \alpha \int_{\Omega \times \Gamma_2}  (h(u_3 + w_3^D(t)) - u_2) v_3 
  dxd\gamma_y 
 \mbox{ for } v_3 \in X. $$
\end{lemma}

The next Lemma is concerned with the continuity property of $\varphi^t$ with respect to $t$. 

\begin{lemma} \label{lem2}
Let $T > 0$. If (A1) $\sim$ (A5) hold, 
$R$ and $Q$ are Lipschitz continuous and bounded on $\R$, 
$\psi$ is Lipschitz continuous on $\R$ 
 and $w_4 \in K(T)$, then for any $r > 0$ there exists 
$b_r \in W^{1,2}(0,T)$ such that for $0 \leq s \leq t \leq T$ and $u \in V$ with $|u|_H \leq r$ it holds  that 
$$ \varphi^t(w_4(t), u) - \varphi^s(w_4(s), u) 
  \leq |b_r(t) - b_r(s)|( 1 + |\varphi^s(w_4(s),u)|). $$
\end{lemma}

Accordingly, by the theory of evolution equations developed  cf. \cite{Ke(Chiba)} and by 
Lemma \ref{lem2}, we can deduce  the solvability of the following problem. 
\begin{lemma} \label{lem3}
Let  $T > 0$. If  (A1) $\sim$ (A5) hold, 
$R$ and $Q$ are Lipschitz continuous and bounded on $\R$, 
$\psi$ is Lipschitz continuous on $\R$ 
 and $w_4 \in K(T)$, then for any $u_0 \in V$ there exists one and only 
$u \in W^{1,2}(0,T;H)$ such that $\varphi^{(\cdot)}(w_4(\cdot), u(\cdot)) \in 
L^{\infty}(0,T)$ and 
\begin{equation}
\frac{d}{dt}u(t) + \partial \varphi^t(w_4(t), u(t)) = f(t) \mbox{ in } H 
  \mbox{ for } t \in [0,T] \mbox{ and } u(0) = u_0, \label{eqe}
\end{equation}
where $f(t) = (0,0, - \partial_t w_3^D(t) +  \nabla d_3 \nabla w_3^D(t))$ 
for  $ t \in [0,T]$. 
\end{lemma}

 Lemma \ref{lem1} and Lemma \ref{lem3} guarantee the existence of a solution of AP($w_4$). 
\begin{lemma} \label{pro1}
Let $T > 0$. If (A1) $\sim$ (A5) hold, 
$R$ and $Q$ are Lipschitz continuous and bounded on $\R$, 
$\psi$ is Lipschitz continuous on $\R$ 
 and $w_4 \in K(T)$, then AP($w_4$) has a unique solution 
$w = (w_1,w_2,w_3)$ on $[0,T]$. 
Precisely speaking, 
$w_1, w_2 \in W^{1,2}(0,T; L^2(\Omega \times Y)) \cap L^{\infty}(0,T; L^2(\Omega; H^1(Y)))$,  
$w_3  \in W^{1,2}(0,T; L^2(\Omega))$, 
$w_3 - w_3^D \in L^\infty(0,T; X)$, (\ref{pd1b}) $\sim$ (\ref{bdb}) hold 
in the usual sense. 
\end{lemma}

{\it Proof. } 
By putting $u_0 = (w_{10},w_{20}, w_{30} - w_3^D(0))$ we see that 
$u_0 \in V$. Then on account of  Lemma \ref{lem3} there exists 
$u = (u_1, u_2, u_3)$ satisfying (\ref{eqe}). Hence, by putting 
$w_1 = u_1$, $w_2 = u_2$ and $w_3 = u_3 + w_3^D$ it is clear that the assertion of this proposition is true. 
\hfill $\Box$

\vskip 12pt
Now, we show a proof of Theorem \ref{Exist} by applying the fixed point argument in case $Q$, $R$ and $\psi$ are Lipschitz continuous. 
\begin{proposition} \label{prro}
Let $T > 0$. If (A1) $\sim$ (A5) hold, 
$R$ and $Q$ are Lipschitz continuous and bounded on $\R$, 
$\psi$ is Lipschitz continuous on $\R$,  
then TP$(R, Q, \psi)$ has a unique solution on $[0,T]$. 

\end{proposition}
{\it Proof. }
Let $T > 0$. 
By Lemma \ref{pro1} for $w_4 \in K(T)$ these exists one and only one 
solution $(w_1, w_2, w_3) \in W^{1,2}(0,T; H)$ of AP($w_4$). 
Then we  can define a mapping $\Lambda_T$ in the following way: 
$\Lambda_T: K(T) \to K(T)$ is given by 
$$ (\Lambda_T w_4)(t) = w_{04} + \int_0^t \eta(w_1(\tau), w_4(\tau)) d\tau 
 \quad \mbox{ for } t \in [0,T]. $$
Obviously, this mapping is well-defined. From now on we shall show 
that $\Lambda_T$ is a contraction mapping for small $T > 0$. 

Let $w_4^{(1)}$, $w_4^{(2)} \in K(T)$, 
$w^{(1)} = (w_1^{(1)},w_2^{(1)},w_3^{(1)})$ 
and $w^{(2)} = (w_1^{(2)},w_2^{(2)},w_3^{(2)})$  
be solutions of AP($w_4^{(1)}$) and AP($w_4^{(2)}$), respectively, 
and put $w_4 = w_4^{(1)} - w_4^{(2)}$, $w = w^{(1)} - w^{(2)} = 
(w_1, w_2, w_3)$.  
It from (\ref{pd1b}) follows that 
$$ \partial_t w_1 - \nabla_y \cdot (d_1 \nabla_y w_1) 
= - \psi(w_1^{(1)} - \gamma w_2^{(1)}) 
+ \psi(w_1^{(2)} - \gamma w_2^{(2)}) \quad \mbox{ in } 
  (0,T ) \times \Omega \times Y. $$
By multiplying it by $w_1$ we have 
\begin{eqnarray}
& & \frac{1}{2} \frac{d}{dt} |w_1|_{L^2(\Omega \times Y)}^2 
  + \int_{\Omega \times Y} d_1|\nabla_y w_1|^2 dxdy \nonumber \\
& &   + \int_{\Omega \times \Gamma_1} 
 (Q(w_4^{(1)}) R(w_1^{(1)}) - Q(w_4^{(2)})R(w_1^{(2)})) w_1 dx d\gamma_y 
 \nonumber  \\
& = & - 
\int_{\Omega \times Y} (\psi(w_1^{(1)} - \gamma w_2^{(1)}) - 
      \psi(w_1^{(2)} - \gamma w_2^{(2)})) w_1 dx dy 
 \quad \mbox{ a.e. on } [0,T]. \label{eqc1}
\end{eqnarray}
Similarly to (\ref{eqc1}), we see that 
\begin{eqnarray*}
& &  \frac{\gamma}{2} \frac{d}{dt} |w_2|_{L^2(\Omega \times Y)}^2 
  + \gamma \int_{\Omega \times Y} d_2|\nabla_y w_2|^2 dxdy \nonumber \\
& = &    \gamma  \alpha \int_{\Omega \times \Gamma_2} 
       (h w_3 - w_2) w_2 dx d\gamma_y 
 \nonumber  \\
&   & + \gamma
  \int_{\Omega \times Y} (\psi(w_1^{(1)} - \gamma w_2^{(1)})
   - \psi(w_1^{(2)} - \gamma w_2^{(2)}) ) w_2 dx dy 
 \quad \mbox{ a.e. on } [0,T],  
\end{eqnarray*}
and 
\begin{eqnarray*}
  \frac{\gamma h}{2} \frac{d}{dt} |w_3|_{L^2(\Omega)}^2 
  + \gamma h  \int_{\Omega} d_3|\nabla w_3|^2 dx \nonumber  
=  - \gamma h \alpha  \int_{\Omega \times \Gamma_2} 
       (h w_3 - w_2) w_3 dx d\gamma_y   \mbox{ a.e. on } [0,T].   
\end{eqnarray*}
Since (A2) implies that $R$ is increasing on $\R$, we obtain 
\begin{eqnarray*}
& &  \int_{\Omega \times \Gamma_1} 
 (Q(w_4^{(1)}) R(w_1^{(1)}) - Q(w_4^{(2)})R(w_1^{(2)})) w_1 dx d\gamma_y 
   \\
& \geq &
  -\int_{\Omega \times \Gamma_1} 
 |Q(w_4^{(1)})  - Q(w_4^{(2)})| |R(w_1^{(2)})| |w_1| dx d\gamma_y
   \quad \mbox{ a.e. on } [0,T]. 
\end{eqnarray*}
Accordingly, it holds that 
\begin{eqnarray*}
& & \frac{1}{2} \frac{d}{dt} |w|_H^2 
 + d_1^0 \int_{\Omega \times Y} |\nabla_y w_1|^2 dxdy
 + \gamma d_2^0 \int_{\Omega \times Y} |\nabla_y w_2|^2 dxdy 
 + \gamma h d_3^0 \int_{\Omega} |\nabla w_3|^2 dx \\
& \leq & 
  \int_{\Omega \times \Gamma_1} 
 |Q(w_4^{(1)})  - Q(w_4^{(2)})| |R(w_1^{(2)})||w_1| dx d\gamma_y  \\
& & -  \int_{\Omega \times Y} (\psi(w_1^{(1)} - \gamma w_2^{(1)})
   - \psi(w_1^{(2)} - \gamma w_2^{(2)}) )  (w_1 - \gamma w_2) dx dy  \\
& \leq & C_Q |R|_{L^{\infty}(\R)} 
\int_{\Omega} |w_1|_{L^2(\Gamma_1)} |w_4|_{L^2(\Gamma_1)} dx \\
& \leq &  C_Y C_Q |R|_{L^{\infty}(\R)} 
\int_{\Omega} (|\nabla_y w_1|_{L^2(Y)} +  |w_1|_{L^2(Y)}) |w_4|_{L^2(\Gamma_1)} dx \\
& \leq &  \frac{d_1^0}{2} \int_{\Omega} |\nabla_y w_1|_{L^2(Y)}^2 dx 
+  \frac{1}{2} \int_{\Omega} |w_1|_{L^2(Y)}^2 dx  
 + C_0 
\int_{\Omega}  |w_4|_{L^2(\Gamma_1)}^2 dx \quad \mbox{ a.e. on } [0,T], 
\end{eqnarray*}
where $C_Q$ is a  Lipschitz constant of $Q$,
$C_0 = (\frac{1}{2d_1^0} + \frac{1}{2}) (C_Y C_Q |R|_{L^{\infty}(\R)})^2$  and 
$C_Y$ is a positive constant such that 
\begin{equation}
|z|_{L^2(\Gamma_1)} \leq C_Y |z|_{H^1(\Omega)} 
 \mbox{ for } z \in H^1(Y). \label{sob}
\end{equation}
 From these inequalities, it follows that 
\begin{eqnarray*}
& & \frac{1}{2} \frac{d}{dt} |w|_H^2 
 + \frac{d_1^0}{2} \int_{\Omega \times Y} |\nabla_y w_1|^2 dxdy
 + \gamma d_2^0 \int_{\Omega \times Y} |\nabla_y w_2|^2 dxdy 
  + \gamma h d_3^0 \int_{\Omega} |\nabla w_3|^2 dxdy \\
& \leq & C_0 \int_{\Omega}  |w_4|_{L^2(\Gamma_1)}^2 dx  
 + \frac{1}{2}  \int_{\Omega} |w_1|_{L^2(Y)}^2 dx  \\
& \leq &  C_1( \int_{\Omega}  |w_4|_{L^2(\Gamma_1)}^2 dx  + |w|^2_H) 
 \quad \mbox{ a.e. on } [0,T], 
\end{eqnarray*}
where $C_1 = C_0 + \frac{1}{2}$. Then by applying Gronwall's inequality we obtain 
\begin{eqnarray} 
& & \frac{1}{2} |w(t)|_H^2 
 +  \int_0^t( \frac{d_1^0}{2} |\nabla_y w_1|_{L^2(\Omega \times Y)}^2 
 + \gamma d_2^0  |\nabla_y w_2|_{L^2(\Omega \times Y)}^2  
   + \gamma h d_3^0   |\nabla w_3|_{L^2(\Omega)}^2) d\tau  
          \nonumber \\
& \leq & 
C_2 \int_0^t |w_4|_{L^2(\Omega \times \Gamma_1)}^2 d\tau \quad \mbox{ for } t \in [0,T],  \label{eq:b10} 
\end{eqnarray}
where $C_2$ is a positive constant. This shows that 
\begin{eqnarray}
& &  |\partial_t(\Lambda w_4^{(1)}) - \partial_t(\Lambda w_4^{(2)})|_{L^2(0,T; L^2(\Omega \times \Gamma_1))} \nonumber \\
&  \leq &  |R(w_1^{(1)}) Q(w_4^{(1)}) - R(w_1^{(2)}) Q(w_4^{(2)})|_{L^2(0,T; L^2(\Omega \times \Gamma_1))} \nonumber \\
& \leq & 
C_{R}|Q|_{L^{\infty}(\R)}|w_1|_{L^2(0,T; L^2(\Omega \times \Gamma_1))}  
+ C_{Q} |R|_{L^{\infty}(\R)}|w_4|_{L^2(0,T; L^2(\Omega \times \Gamma_1))}) 
\nonumber \\
& \leq &  (1 + C_Y) C_{R}|Q|_{L^{\infty}(\R)} 
      |w_1|_{L^2(0,T; L^2(\Omega; H^1(Y)))}  
+ C_{Q} |R|_{L^{\infty}(\R)} |w_4|_{L^2(0,T; L^2(\Omega \times \Gamma_1))}) 
\nonumber \\
& \leq & C_3 |w_4|_{L^2(0,T; L^2(\Omega \times \Gamma_1))} \nonumber \\
& \leq & C_3 |\int_0^t \partial_{\tau} w_{4} d\tau |_{L^2(0,T; L^2(\Omega \times \Gamma_1))} \nonumber \\
& \leq & 
 C_3 T |\partial_t w_{4}|_{L^2(0,T; L^2(\Omega \times \Gamma_1))}, \label{eq:73}\end{eqnarray}
where $C_R$ is a Lipschitz constant of $R$ and $C_3$ is a positive constant depending only on  $C_2$. 
Moreover,  by using (\ref{eq:73}) we have 
\begin{eqnarray*}
 |\Lambda w_4^{(1)} - \Lambda w_4^{(2)}|_{L^2(0,T; L^2(\Omega \times \Gamma_1))}& \leq & |\int_0^t \partial_t(\Lambda w_4^{(1)} - \Lambda w_4^{(2)}) d\tau|_{L^2(0,T; L^2(\Omega \times \Gamma_1))}  \\
& \leq & T |\partial_t(\Lambda w_4^{(1)}) - \partial_t(\Lambda w_4^{(2)})|_{L^2(0,T; L^2(\Omega \times \Gamma_1))} \\
& \leq & 
C_4 T^2 |\partial_t w_{4}|_{L^2(0,T; L^2(\Omega \times \Gamma_1))},
\end{eqnarray*}
where $C_4$ is a positive constant.

Hence, if $T_1$ is sufficiently small, then 
$\Lambda_{T_1}$ is a contraction mapping. 
Namely, Banach's fixed point theorem shows that TP$(R, Q, \psi)$ has a solution on $[0,T_1]$. Furthermore, the choice of $T_1$ is independent of initial values so that we conclude that TP$(R, Q, \psi)$ has a solution on $[0,T]$. 

The uniqueness of a solution is a direct consequence of the Banach's fixed point theorem. 
\hfill $\Box$

\vskip 12pt
\begin{lemma} \label{bdd}
Let $T > 0$ and assume that  (A1) $\sim$ (A5) hold, 
$R$ and $Q$ are Lipschitz continuous and bounded on $\R$, 
$\psi$ is Lipschitz continuous on $\R$. 
If $(w_1, w_2, w_3, w_4)$ is a solution of TP$(R, Q, \psi)$ on $[0,T]$,
then it holds that 
$$ 0 \leq w_i \leq M_i \mbox{ a.e. on } (0,T) \times \Omega \times Y \mbox{ for  } i= 1,2,$$ 
$$ 0 \leq w_3 \leq M_3 \mbox{ a.e. on } (0,T) \times \Omega, \mbox{ and  } 
0 \leq w_4 \leq M_4 \mbox{ a.e. on } (0,T) \times \Omega \times \Gamma_1, 
$$ 
where $M_i$ is the positive constant defined in Theorem \ref{Exist} for each 
$ i= 1, 2, 3, 4$.  
\end{lemma}
The proof of this lemma is quite similar to that of \cite[Theorem 4.4]{two-scale1} so that we omit a proof of Lemma \ref{bdd}. 

\vskip 12pt
By using Proposition \ref{prro} and Lemma \ref{bdd}, we can now prove Theorem \ref{Uniqueness} and Theorem \ref{Exist}. 

{\it Sketch of the proofs of Theorem \ref{Uniqueness} and Theorem \ref{Exist}. }
First, for $m > 0$ we put
$$ R_m(r) = \left\{ \begin{array}{cl} 
              R(m) & \mbox{ for } r > m, \\
               R(r) & \mbox{ otherwise, } 
             \end{array} \right.
Q_m(r) = \left\{ \begin{array}{cl} 
Q(m) & \mbox{ for } r > m, \\
Q(r) & \mbox{ for } 0 \leq r \leq m, \\
               Q(0) & \mbox{ otherwise, } 
             \end{array} \right.
 $$
and 
$$
\psi_m(r) = \left\{ \begin{array}{cl} 
\psi(m) & \mbox{ for } r > m, \\
\psi(r) & \mbox{ for } 0 \leq r \leq m, \\
               \psi(-m) & \mbox{ otherwise. } 
             \end{array} \right.
 $$
Then, Proposition \ref{prro} implies that TP$(R_m, Q_m, \psi_m)$ has a solution $(w_{1m}, w_{2m}, w_{3m}, w_{4m})$ on $[0,T]$ for any $T > 0$. Moreover, 
on account of Lemma \ref{bdd} we can choose a positive constant $m > 0$ such that $R_m(w_{1m}) = R(w_{1m})$, $Q_m(w_{1m}) = Q(w_{1m})$ and 
$\psi_m(w_{1m} - \gamma w_{2m})  = \psi(w_{1m} - \gamma w_{2m})$, 
since $M_i$ is independent of $m$ for each $i = 1,2,3,4$. 
This shows the conclusion of Theorem \ref{Exist}.  
Moreover, the uniqueness of a solution of TP$(R, Q, \psi)$ is a direct consequence of Proposition \ref{prro}. 
\hfill $\Box$

\section{The quest of the large-time behavior} \label{LARGE}

To prove Theorem \ref{large}, throughout this section 
we always assume  (A1) $\sim$ (A5) and (A6), and 
denote a solution of TP$(R, Q, \psi)$ on $[0,\infty)$ 
by $(w_1, w_2, w_3, w_4)$.  
Clearly, by putting $f := \partial_t w_3^D - \nabla d_3 \nabla w_3^D$ 
there exists $f_{\infty} \in L^{\infty}(\Omega)$ such that 
$f(t) \to f_{\infty}$ in $L^1(\Omega)$ as $t \to \infty$. 
Also, let $m_0 = \max\{M_1,  M_2, M_3, M_4\}$. Then for some $m_1 \geq m_0$ 
it holds that 
$R_{m_1}(z_1) = R(z_1)$, $\psi_{m_1}(z_1 - \gamma z_2) = \psi(z_1 - \gamma z_2)$ for $z = (z_1, z_2, z_3) \in U(m_0)$, where 
$U(m_0) = \{z = (z_1, z_2, z_3) \in H| 0 \leq z_1 \leq m_0, 0 \leq z_2 \leq m_0 \mbox{ a.e. on } \Omega \times Y, 0 \leq z_3 \leq m_0 \mbox{ a.e. on } \Omega\}$. Moreover, for simplicity, we write $R$, $Q$ and $\psi$ as $R_{m_1}$, $Q_{m_1}$ and $\psi_{m_1}$, respectively, in this section, and set 
$$ \ell_1 = R(M_1), \ell_2 = |Q(w_4)|_{L^{\infty}((0,\infty) \times \Omega \times \Gamma_1)}, \ell_3 = \hat{R}(M_1), 
\ell_4 = \sup \{ |Q'(r)|: |r| \leq m_0|\}. $$

First, we show the convergence of $w_4(t)$ as $t \to \infty$. 
\begin{lemma} \label{w4} 
$w_4(t)$ converges to some function $w_{4\infty} \in L^{\infty}(\Omega \times \Gamma_1)$ in $L^1(\Omega \times \Gamma_1)$ as $t \to \infty$, $0 \leq w_{4\infty} \leq m_0$ a.e. on $\Omega \times \Gamma_1$, and 
\begin{equation}
 \partial_t w_4 \in L^1(0,\infty; L^1(\Omega \times \Gamma_1)) \cap 
L^2(0,\infty; L^2(\Omega \times \Gamma_1)). \label{tb}
\end{equation} 
\end{lemma}
{\it Proof. } By (A2) it is obvious that $\partial_t w_4 = \eta(w_1, w_4) \geq 0$  a.e. on $(0,\infty) \times \Omega \times \Gamma_1$. Then Lemma \ref{bdd} and the Lebesgue monotone convergence theorem imply that 
$w_4(t) \to w_{4\infty}$ in $L^1(\Omega \times \Gamma_1)$ as $t \to \infty$, 
where $w_{4\infty} \in L^{\infty}(\Omega \times \Gamma_1)$.  

It is easy to see that 
\begin{eqnarray*}
\int_0^T \int_{\Omega \times \Gamma_1} 
 |\partial_t w_4| dx d\gamma_y dt 
& = & 
\int_0^T \int_{\Omega \times \Gamma_1} 
 \partial_t w_4 dx d\gamma_y dt  \\
& = & 
\int_{\Omega \times \Gamma_1} 
  (w_4(T) - w_{04})  dx d\gamma_y \quad \mbox{ for } T  > 0, 
\end{eqnarray*}
and 
\begin{eqnarray*}
\int_0^T \int_{\Omega \times \Gamma_1} 
 |\partial_t w_4|^2 dx d\gamma_y  dt 
& \leq  &  \ell_1 \ell_2  \int_0^T \int_{\Omega \times \Gamma_1} 
 |\partial_t w_4| dx d\gamma_y dt  
 \mbox{ for } T  > 0.  
\end{eqnarray*}
Then these estimates lead to (\ref{tb}). \hfill $\Box$

\vskip 12pt
Now, we provide uniform estimates on the time derivative of  solutions.
\begin{lemma} \label{es2}
It holds that 
$$ \partial_t w_1, \partial_t w_2 \in L^2(0, \infty; L^2(\Omega \times Y)), 
\partial_t \tilde{w}_3 \in L^2(0, \infty; L^2(\Omega)), $$
$$ w_1, w_2 \in L^{\infty}(0, \infty; L^2(\Omega; H^1(Y))), 
 \tilde{w}_3 \in L^{\infty}(0, \infty; X), $$
where $\tilde{w}_3 = w_3 - w_{3}^D$.
\end{lemma} 
{\it Proof. }  
By multiplying (\ref{pd1}) by $\partial_t w_1$,  we can obtain 
\begin{eqnarray*}
& & \int_{\Omega \times Y}  |\partial_t w_1|^2 dx dy + 
\frac{1}{2} \frac{d}{dt} \int_{\Omega \times Y} d_1 |\nabla_y w_1|^2 dxdy 
+ \frac{d}{dt} \int_{\Omega \times \Gamma_1} 
  Q(w_4) \hat{R}(w_1) dx d\gamma_y \\
& = & - \int_{\Omega \times Y} \psi(w_1 - \gamma w_2) \partial_t w_1 dxdy 
+ 
\int_{\Omega \times \Gamma_1} 
           Q'(w_4) \partial_t w_4 \hat{R}(w_1) dx d\gamma_y 
 \quad \mbox{ a.e. on } [0,\infty). 
\end{eqnarray*}
Next, we multiply (\ref{pd2}) by $\gamma \partial_t w_2$. Then we see that 
\begin{eqnarray*}
& &  \gamma \int_{\Omega \times Y}  |\partial_t w_2|^2 dx dy + 
\frac{\gamma}{2} \frac{d}{dt} \int_{\Omega \times Y} d_2|\nabla_y w_2|^2 dxdy 
 \\
& & 
+ \gamma \alpha  \int_{\Omega \times \Gamma_2} 
  (w_2 - h(\tilde{w}_3 + w_3^D) ) \partial_t w_2 dx d\gamma_y \\
& = & 
 \gamma  \int_{\Omega \times Y} \psi(w_1 - \gamma w_2) \partial_t w_2 dxdy 
 \quad \mbox{ a.e. on } [0,\infty). 
\end{eqnarray*}
Easily, we have 
$$ 
\partial_t \tilde{w}_3 - \nabla d_3 \nabla \tilde{w}_3 = 
 -\alpha\int\limits_{\Gamma_2}\big( h (\tilde{w}_3 + w_3^D) - w_2\big)d\gamma_y
- f  \quad \mbox{  a.e. on  } (0,\infty) \times \Omega, $$
where $f = \partial_t w_3^D - \nabla d_3\nabla w_3^D$. 
Here, we multiply it by $\gamma h \partial_t \tilde{w}_3$ and observe that 
\begin{eqnarray*}
& & \gamma h \int_{\Omega} |\partial_t \tilde{w}_3|^2 dx 
 + \frac{\gamma h}{2} \frac{d}{dt} \int_{\Omega} d_3 |\nabla \tilde{w}_3|^2 dx \\
& = & - \alpha \gamma h \int_{\Omega \times \Gamma_2} 
  (h (\tilde{w}_3 + w_3^D) - w_2) \partial_t \tilde{w}_3 dx d\gamma_y 
- \gamma h \int_{\Omega} f \partial_t \tilde{w}_3 dx \\
& = & - \alpha \gamma h \int_{\Omega \times \Gamma_2} 
  (h (\tilde{w}_3 + w_3^D) - w_2) (\partial_t \tilde{w}_3 + \partial_t w_3^D) 
 dx d\gamma_y \\
& & + \alpha \gamma h \int_{\Omega \times \Gamma_2} 
  (h (\tilde{w}_3 + w_3^D) - w_2)  \partial_t w_3^D  dx d\gamma_y 
- \gamma h \frac{d}{dt} \int_{\Omega} f \tilde{w}_3 dx 
+ \gamma h \int_{\Omega} \partial_t f \tilde{w}_3 dx
\end{eqnarray*}
a.e. on $[0,\infty)$. By adding these equations we obtain 
\begin{eqnarray}
 & & \int_{\Omega \times Y}  |\partial_t w_1|^2 dx dy 
 + \gamma \int_{\Omega \times Y}  |\partial_t w_2|^2 dx dy 
 + \gamma h \int_{\Omega}  |\partial_t \tilde{w}_3|^2 dx \nonumber \\
& &  + \frac{1}{2} \frac{d}{dt} 
\left( \int_{\Omega \times Y} d_1 |\nabla_y w_1|^2 dxdy 
+ \gamma \int_{\Omega \times Y} d_2 |\nabla_y w_3|^2 dxdy 
+ \gamma h \int_{\Omega}  d_3 |\nabla  \tilde{w}_3|^2 dx \right) \nonumber \\
& & + \frac{d}{dt} \int_{\Omega \times \Gamma_1} 
  Q(w_4) \hat{R}(w_1) dx d\gamma_y 
+  \gamma \alpha \frac{d}{dt}  \int_{\Omega \times \Gamma_2} 
  |w_2 - h(\tilde{w}_3 + w_3^D)|^2  dx d\gamma_y \nonumber \\
& &  +  \frac{d}{dt} \int_{\Omega \times Y
} \hat{\psi}(w_1 - \gamma w_2)  dxdy  \nonumber \\
& = & 
\int_{\Omega \times \Gamma_1} 
           Q'(w_4) \partial_t w_4 \hat{R}(w_1) dx d\gamma_y 
+ \alpha \gamma h \int_{\Omega \times \Gamma_2}
  (h(\tilde{w}_3 + w_3^D) - w_2) \partial w_{3}^D dx d\gamma_y 
  \label{eq13} \\
& & - \gamma h \frac{d}{dt} \int_{\Omega} f \tilde{w}_3 dx 
+ \gamma h \int_{\Omega} \partial_t f \tilde{w}_3 dx \ \ 
(=: \sum_{i=1}^4 I_i(\cdot) ) 
 \quad \mbox{ a.e. on } [0,\infty).  \nonumber
\end{eqnarray}
Since $$ I_1(t) \leq 
 \ell_3 \ell_4 
 \int_{\Omega \times \Gamma_1} |\partial_t w_4(t)| dxd\gamma_y  
 \mbox{ for a.e. } t \geq 0, $$
it holds that $I_1 \in L^1(0,\infty)$. 
It is clear that 
\begin{eqnarray*}
I_2 & = & 
  \alpha \gamma h \int_{\Omega \times\Gamma_2} (h w_3 - w_2) \partial_t w_3^D 
    dx d\gamma_y \\
& \leq & \alpha \gamma h (h M_3 + M_2) 
      \int_{\Omega \times\Gamma_2} |\partial_t w_3^D| dx d\gamma_y 
 \quad \mbox{ a.e. on } [0,\infty)  
\end{eqnarray*}
so that $I_2 \in L^1(0,\infty)$. Also, we see that 
$$ \int_0^t I_3(\tau) d\tau 
  =   \gamma h( - \int_{\Omega} f(t) \tilde{w}_3(t) dx 
+  \int_{\Omega} f(0) (w_{30} + w_{3}^D(0)) dx ) \leq B_1 \mbox{ for } t \geq 0,$$
where $B_1$ is a positive constant, and 
$$ I_4(t) \leq \gamma h (M_3 + |w_3^D|_{L^{\infty}((0,\infty) \times \Omega)}) 
 \int_{\Omega} |\partial_t f(t)| dx \quad 
 \mbox{ for } t \geq 0 
$$
so that $I_4 \in L^{\infty}(0,\infty)$. 
Hence, we have proved this lemma. \hfill $\Box$ 

\vskip 12pt 
To describe the large-time behavior of the solution, we introduce the following notations: We put 
$w(t) := (w_1(t), w_2(t), w_3(t)) \in H$ for $t \geq 0$, 
$$ \omega(w_0) =\left\{ z \in H \left|  \begin{array}{l} 
w(t_n) \to z \mbox{ weakly in } H \mbox{ as } 
 n \to \infty 
\mbox{ for some sequence } \{t_n\}  \\
\mbox{with } t_n \to \infty \mbox{ as } n \to \infty \end{array}
\right. \right\},  $$
where $w_0 = (w_{10}, w_{20}, w_{30})$, 
$$ \varphi_1^t(u) = \varphi^t(w_4(t), u) - \gamma h \int_{\Omega} f(t) u_3 dx 
 \mbox{ for } u = (u_1, u_2, u_3) \in H, $$ 
\begin{eqnarray*}
& & \varphi_1^{\infty}(u) \\
& = & 
\left\{ \begin{array}{l} \displaystyle
\frac{1}{2} \int_{\Omega \times Y} d_1 |\nabla_y u_1|^2 dxdy 
+  \int_{\Omega \times \Gamma_1} Q(w_{4\infty})  \hat{R}(u_1) dxd\gamma_y 
+ \frac{\gamma}{2} \int_{\Omega \times Y} d_2 |\nabla_y u_2|^2 dxdy \\[0.3cm]
\displaystyle
+  \int_{\Omega \times Y} \hat{\psi} (u_1 - \gamma u_2) dxdy  
 + \frac{\gamma \alpha}{2} \int_{\Omega \times \Gamma_2} 
        |h(u_3 + w_{3\infty}^D) - u_2|^2 dxd\gamma_y 
 \\[0.3cm]
\displaystyle
  + \frac{\gamma h}{2} \int_{\Omega} d_3|\nabla u_3|^2 dx 
 - \gamma h \int_{\Omega} f_{\infty} u_3 dx \quad 
\mbox{ if } u = (u_1, u_2, u_3) \in V, \\
\infty  \quad \quad  \mbox{ otherwise, }
\end{array} \right. 
\end{eqnarray*}
and 
$$ F(\varphi_1^{\infty}) 
  = \{ z \in H| \varphi_1^{\infty}(z) = 
 \min_{u \in H} \varphi_1^{\infty}(u)\}. $$
Clearly, the similar results to Lemma \ref{lem1} hold for $\varphi_1^{\infty}$. Here, we note that $w$ satisfies (S2), (S3) and (S4) if and only if 
$\tilde{w} = (w_1, w_2, \tilde{w}_3)$ with $\tilde{w}_3 = w_3 - w_3^D$ is a solution of the following evolution equation: 
\begin{equation}
\tilde{w}_t(t) + \partial \varphi_1^t(\tilde{w}(t)) = 0 \mbox{ in } H \mbox{ for a.e. } t \geq 0.  
\label{eq101}
\end{equation}

Moreover, since by Theorem \ref{Exist} $\{w(t)\}_{t\geq 0}$ is bounded in $H$, 
there exist a sequence $\{t_n\}$ with $t_n \to \infty$ as $n \to \infty$ 
and $w_{\infty} = (w_{1\infty}, w_{2\infty}, w_{3\infty}) \in U(m_0)$ such that 
$$ w(t_n) \to w_{\infty} \mbox{ weakly in } H \mbox{ as } n \to \infty, $$
that is, $w_{\infty} \in \omega(w_0)$. 

The next Lemma guarantees the existence of a solution to 
SP$(w_{4\infty}, w_{3\infty}^D$).  
\begin{lemma} \label{spp}
$w_{\infty} \in F(\varphi_1^{\infty})$ and 
$w_{\infty}$ is a solution of SP$(w_{4\infty}, w_{3\infty}^D$).  
\end{lemma}
{\it Proof. } 
By integrating (\ref{eq13}) over $[s, t]$ with $0 \leq s \leq t$ 
we have 
\begin{eqnarray*}
 \varphi_1^t(\tilde{w}(t)) - \varphi_1^s(\tilde{w}(s)) 
+ \int_s^t |\partial_t \tilde{w}|_{H}^2  d\tau  
 \leq  \sum_{i = 1}^4 \int_s^t I_i d\tau \quad \mbox{ for } 
 0 \leq s \leq t,  
 \end{eqnarray*}
where $\tilde{w} = (w_1, w_2, \tilde{w}_3)$. 
It is obvious that  the function $t \to \varphi_1^t(\tilde{w}(t))
+ \int_0^t |\partial_t \tilde{w}|_{H}^2  d\tau 
- \sum_{i = 1}^4 \int_0^t I_i d\tau$ is non-increasing on $[0,\infty)$ and 
$\varphi_1^t(\tilde{w}(t))) \geq 0$ for $t \geq 0$. 
Then  by  Lemma \ref{es2} we see that 
$\lim_{t \to \infty} \varphi_1^t(\tilde{w}(t))$ exists so that we can put 
$q_0 = \lim_{t \to \infty} \varphi_1^t(\tilde{w}(t))) \geq 0$. 

Next, we show that 
\begin{equation}
\varphi_1^{\infty}(\tilde{w}_{\infty}) \leq 
 \liminf_{n \to \infty} \varphi_1^{t_n}(\tilde{w}(t_n)), 
\label{eq41}
\end{equation}
where $\tilde{w}_{\infty} =  (w_{1\infty}, w_{2\infty}, 
w_{3\infty} - w_{3\infty}^D)$. 
In fact, we observe that for each $n$ 
\begin{eqnarray*}
& & |\varphi_1^{\infty}(\tilde{w}(t_n))
- \varphi_1^{t_n}(\tilde{w}(t_n))| \\
& \leq & 
\int_{\Omega \times \Gamma_1} |Q(w_{4\infty}) \hat{R}(w_1(t_n))
    -  Q(w_{4}(t_n)) \hat{R}(w_1(t_n)) | dxd\gamma_y \\
& & + \frac{\gamma \alpha}{2} \int_{\Omega \times \Gamma_2}
  | |h(\tilde{w}_3(t_n) + w_{3\infty}^D) - w_2(t_n)|^2
 - |h(\tilde{w}_3(t_n) + w_{3}^D(t_n)) - w_2(t_n))|^2|  dx d\gamma_y \\
& & + h \int_{\Omega} |(f_{\infty} - f(t_n)) \tilde{w}_3(t_n)| dx \\
& \leq & \ell_3 \ell_4
\int_{\Omega \times \Gamma_1} |w_{4\infty} - w_{4}(t_n)| dxd\gamma_y 
+  h \int_{\Omega} |f_{\infty} - f(t_n)||\tilde{w}_3(t_n)| dx \\
& & + \frac{\gamma \alpha}{2} \int_{\Omega \times \Gamma_2}
  |h(w_{3\infty}^D - w_{3}^D(t_n))
 (2h \tilde{w}_3(t_n) + h( w_{3}^D(t_n) + w_{3\infty}^D) - 2w_2(t_n)) |
   dx d\gamma_y. \\
\end{eqnarray*}
Then, for some positive constant $C_5$, we obtain: 
\begin{eqnarray*}
& & |\varphi_1^{\infty}(\tilde{w}(t_n))
- \varphi_1^{t_n}(\tilde{w}(t_n))| \\
& \leq & C_5 \int_{t_n}^{\infty} 
(\int_{\Omega \times \Gamma_1} 
       |\partial_t w_4| dxd\gamma_y 
 + \int_{\Omega \times \Gamma_2} 
       |\partial_t w_3^D| dxd\gamma_y 
 + \int_{\Omega} 
       |\partial_t f| dx )dt \mbox{ for } n
\end{eqnarray*}
so that 
\begin{equation}
\lim_{n \to \infty} (\varphi_1^{\infty}(\tilde{w}(t_n))
- \varphi_1^{t_n}(\tilde{w}(t_n))) = 0. 
\label{con1}
\end{equation}
Hence,  we see that 
\begin{eqnarray*}
 \liminf_{n \to \infty} 
         \varphi_1^{t_n}(\tilde{w}(t_n))
& \geq & \liminf_{n \to \infty} 
 ( \varphi_1^{t_n}(\tilde{w}(t_n)) - 
           \varphi_1^{\infty}(\tilde{w}(t_n))
+ \liminf_{n \to \infty} 
         \varphi_1^{\infty}(\tilde{w}(t_n)) 
 \geq  \varphi_1^{\infty}(\tilde{w}_{\infty}). 
\end{eqnarray*}
Thus (\ref{eq41}) is true. Moreover, it is clear that 
$\varphi_1^{\infty}(\tilde{w}_{\infty}) \leq q_0$ so that 
$\tilde{w}_{\infty} \in V$.

Based on Lemma \ref{es2}, we can take a subsequence $\{t_n'\}$ with 
$t_n' \to \infty$ as $n \to \infty$ such that 
\begin{equation}
 n \leq t_n' \leq n + 1 \mbox{ for each } n \mbox{ and } 
   \partial_t \tilde{w}(t_n') \to 0 \mbox{ in } H \mbox{ as } n \to \infty. 
\label{tn}
\end{equation}
Let $z \in V$. Similarly to (\ref{con1}), we can prove that
$\varphi_1^{t_n'}(z) \to \varphi_1^{\infty}(z)$ 
as $n \to \infty$. Also, by  (\ref{eq101})  we have 
\begin{eqnarray*}
 (-\tilde{w}_t(t_n'), z - w(t_n'))_H 
& = & ( \partial \varphi_1^{t_n'}(w(t_n')), z - \tilde{w}(t_n'))_H \\
& \leq & \varphi_1^{t_n'}(z) - \varphi_1^{t_n'}(\tilde{w}(t_n')) 
 \quad \mbox{ for each } n.
\end{eqnarray*}
By letting $n \to \infty$ in the above inequality,  we obtain 
$0 \leq \varphi_1^{\infty}(z) - q_0$ for any $z \in V$. 
This shows that $q_0 = \min \varphi^{\infty}$, 
that is,  and $\tilde{w}_{\infty} \in F(\varphi^{\infty}_1)$. 
Moreover, we see that 
$0 = \partial \varphi^{\infty}(\tilde{w}_{\infty})$. 
Hence, Lemma \ref{lem1} together with $w_{\infty} \in U(m_0)$ 
implies the conclusion of this Lemma. 
\hfill $\Box$

\vskip 12pt
The next Lemma guarantees the uniqueness of  solutions to 
SP$(w_{4\infty}, w_{3\infty}^D$). 
\begin{lemma} \label{SPuni}
If $(\psi(r) - \psi(r'))(r - r') \geq \mu |r - r'|^{p+1}$ for $r, r' \in \R$, where $\mu > 0$ and $p \geq 1$, then 
SP($w_{4\infty}$, $w_{3\infty}^D$) has at most one solution.  
\end{lemma}

{\it Proof. } 
Let ${w}_{\infty}^{(1)} (:= (w_1^{(1)}, w_2^{(1)}, {w}_3^{(1)}))$ 
and ${w}_{\infty}^{(2)} (:= (w_1^{(2)}, w_2^{(2)}, {w}_3^{(2)}))$ 
be solutions of  SP($w_{4\infty}$, $w_{3\infty}^D$) and 
${w}_{\infty}^{(1)} - {w}_{\infty}^{(2)} = 
(w_{1\infty}, w_{2\infty}, {w}_{3\infty})$.

By Definition \ref{def_SP} we see that 
\begin{eqnarray*}
0 & = & 
\int_{\Omega \times Y} d_1 |\nabla_y w_{1\infty}|^2  dxdy 
  + \int_{\Omega \times \Gamma_1} Q(w_{4\infty}) 
  (R(w_1^{(1)}) - R(w_1^{(2)}) ) w_{1\infty} dxd\gamma_y
\\
& & + \int_{\Omega \times Y}
 ( \psi(w_1^{(1)} - \gamma w_2^{(1)}) -  \psi(w_1^{(2)} - \gamma w_2^{(2)}) )
  (w_{1\infty} - \gamma w_{2\infty})  dxdy \\
&  & +  \gamma \int_{\Omega \times Y} d_2 |\nabla_y w_{2\infty}|^2 dxdy 
  + \alpha \gamma 
 \int_{\Omega \times \Gamma_2} |h {w}_{3\infty}  - w_{2\infty}|^2
   dxd\gamma_y 
 + \alpha \gamma \int_{\Omega} d_3 |\nabla {w}_{3\infty}|^2 dx.
\end{eqnarray*}
Using of the monotonicity of $R$ and $\psi$, 
we obtain $\nabla {w}_{3\infty} = 0$ a.e. on $\Omega$. Since 
${w}_{3\infty} = 0$ a.e.  on $\Gamma_D$, we obtain 
${w}_{3\infty} = 0$ a.e. on $\Omega$. 
Accordingly, we have $w_2 = 0$ a.e.  on $\Omega \times \Gamma_2$. 
Immediately, we infer that 
$w_2 = 0$ a.e. on $\Omega \times Y$. Moreover, it follows that 
$ \int_{\Omega \times Y} d_1 |\nabla_y w_{1\infty}|^2  dxdy + 
 \mu \int_{\Omega \times Y} |w_{1\infty} - \gamma w_{2\infty}|^{p+1} dxdy = 
0$. This implies the uniqueness of a solution of the stationary problem. 
\hfill $\Box$ 

\vskip 12pt
Now, we accomplish the proof of Theorem \ref{large}. 

{\it Proof of Theorem \ref{large}. } 
\newline
Considering the above arguments, it is sufficient to show (\ref{as3}). Let 
$\{t_n'\}$ be a sequence satisfying (\ref{tn}). 
First, we observe that 
\begin{eqnarray*}
& & 
(\partial \varphi_1^{t_n'}(\tilde{w}(t_n'))  -
        \partial \varphi_1^{\infty}(\tilde{w}_{\infty}), 
    \tilde{w}(t_n') -  \tilde{w}_{\infty})_H \\
& \geq & 
\int_{\Omega \times Y} d_1 |\nabla_y (w_1(t_n') - w_{1\infty})|^2  dxdy \\
& & 
  + \int_{\Omega \times \Gamma_1} 
 (Q(w_{4}(t_n')) R(w_1(t_n')) - Q(w_{4\infty}) R(w_{1\infty}) ) 
    (w_1(t_n') - w_{1\infty}) dxd\gamma_y
\\
& & + \mu \int_{\Omega \times Y}
 |(w_1(t_n') - w_{1\infty}) - \gamma (w_2(t_n') -  w_{2\infty})|^{p+1}  dxdy \\
&  & +  \gamma \int_{\Omega \times Y} d_2 
 |\nabla_y (w_2(t_n') - w_{2\infty})|^2 dxdy \\
& &   + \alpha \gamma 
 \int_{\Omega \times \Gamma_2} |h (\tilde{w}_3(t_n') - \tilde{w}_{3\infty})  - 
     (w_2(t_n') - w_{2\infty})|^2    dxd\gamma_y \\
& &   + \alpha \gamma h
 \int_{\Omega \times \Gamma_2}
 ({w}_3^D(t_n') - {w}_{3\infty}^D) \{ 
     (w_2(t_n') - w_{2\infty}) - h(\tilde{w}_3(t_n') - \tilde{w}_{3\infty})\} 
    dxd\gamma_y \\
& &  + \gamma h \int_{\Omega} d_3 
  |\nabla (\tilde{w}_3(t_n') - \tilde{w}_{3\infty})|^2 dx
+ \gamma h \int_{\Omega} (f(t_n') - f_{\infty})(\tilde{w}_3(t_n') - \tilde{w}_{3\infty}) dx 
\\
& = : &  \sum_{i= 1}^8 J_{in} \quad \mbox{ for each } n. 
\end{eqnarray*}
For each $n$,   it is easy to see that 
\begin{eqnarray*}
& & J_{1n} + J_{4n} + J_{7n} \\
& \geq & 
d_1^0 \int_{\Omega \times Y}  |\nabla_y (w_1(t_n') - w_{\infty})|^2dxdy 
+  \gamma d_2^0 \int_{\Omega \times Y} 
    |\nabla_y (w_2(t_n') - w_{2\infty})|^2 dxdy \\
& & +  \gamma h d_3^0 \int_{\Omega} 
  |\nabla (\tilde{w}_3(t_n') - \tilde{w}_{3\infty})|^2 dx, 
\end{eqnarray*}
and 
\begin{eqnarray*}
J_{2n}  
& \geq &  \int_{\Omega \times \Gamma_1} 
 (Q(w_{4}(t_n'))  - Q(w_{4\infty})) R(w_{1\infty})  
    (w_1(t_n') - w_{1\infty}) dxd\gamma_y.   
\end{eqnarray*}
On the other hand, we have 
\begin{equation}
 (\partial \varphi_1^{t_n'}(\tilde{w}(t_n'))  -
        \partial \varphi_1^{\infty}(\tilde{w}_{\infty}), 
    \tilde{w}(t_n') -  \tilde{w}_{\infty})_H =
  (- \tilde{w}_t(t_n'), \tilde{w}(t_n') -  \tilde{w}_{\infty})_H 
   \quad \mbox{ for each } n.  \label{eq53}
\end{equation}
Then we obtain 
\begin{eqnarray}
& & d_1^0 \int_{\Omega \times Y}  |\nabla_y (w_1(t_n') - w_{\infty})|^2dxdy 
+  \gamma d_2^0 \int_{\Omega \times Y} 
    |\nabla_y (w_2(t_n') - w_{2\infty})|^2 dxdy  \nonumber \\
& & + \alpha \gamma d_3^0 \int_{\Omega} 
  |\nabla (\tilde{w}_3(t_n') - \tilde{w}_{3\infty})|^2 dx \nonumber \\
& & + \mu \int_{\Omega \times Y}
  |(w_1(t_n') - w_{1\infty}) - \gamma (w_2(t_n') -  w_{2\infty})|^{p+1} 
  dxdy \nonumber \\
& &   + \alpha \gamma 
 \int_{\Omega \times \Gamma_2} |h (\tilde{w}_3(t_n') - \tilde{w}_{3\infty})  - 
     (w_2(t_n') - w_{2\infty})|^2    dxd\gamma_y \nonumber \\
& \leq & E(t_n') +       |\tilde{w}_t(t_n')|_H |\tilde{w}(t_n') -  \tilde{w}_{\infty}|_H 
    \quad \mbox{ for each } n,  \label{eq61}
\end{eqnarray}
where 
\begin{eqnarray*}
E(t) &=&  \int_{\Omega \times \Gamma_1} 
 |Q(w_{4}(t))  - Q(w_{4\infty})| |R(w_{1\infty})|  
    |w_1(t) - w_{1\infty}| dxd\gamma_y  \\
& &  + \alpha \gamma h  \int_{\Omega \times \Gamma_2}
 ({w}_3^D(t_n') - {w}_{3\infty}^D) \{ 
     (w_2(t_n') - w_{2\infty}) - h(\tilde{w}_3(t_n') - \tilde{w}_{3\infty})\} 
    dxd\gamma_y  \\
& &  + 
\gamma h \int_{\Omega} |f(t) - f_{\infty}||\tilde{w}_3(t) - \tilde{w}_{3\infty}|dx  \quad \mbox{  for } t \geq 0. 
\end{eqnarray*}
By letting $n \to \infty$ in the above inequality, we  infer that 
$\tilde{w}_3(t_n') \to \tilde{w}_{3\infty}$ in $X$, 
$w_2(t_n') \to w_{2\infty}$ and $w_1(t_n') \to w_{1\infty}$ in 
$L^2(\Omega; H^1(Y))$ as $n \to \infty$. In particular, 
$\tilde{w}(t_n') \to  \tilde{w}_{\infty}$ in $H$ as $n \to \infty$.  

From (\ref{eq53}),  it follows that 
$$
- (\partial \varphi_1^{t}(\tilde{w}(t))  -
        \partial \varphi_1^{\infty}(\tilde{w}_{\infty}), 
    \tilde{w}(t) -  \tilde{w}_{\infty})_H = \frac{1}{2} 
 \frac{d}{dt} |\tilde{w}(t) -  \tilde{w}_{\infty}|_H^2 
   \quad \mbox{ for a.e. } t \geq 0.  
$$
Then, similarly to (\ref{eq61}),  we can show that 
$$
\frac{1}{2} 
     \frac{d}{dt} |\tilde{w}(t) -  \tilde{w}_{\infty}|_H^2 
 \leq E(t)  
    \quad \mbox{ for a.e.  } t \geq 0.   $$
By integrating it over $[t_n',t]$ with $t_n' \leq t \leq t_n'+2$ 
for $n$, we have 
\begin{equation}
  \frac{1}{2} |\tilde{w}(t) -  \tilde{w}_{\infty}|_H^2  
 \leq   \frac{1}{2} |\tilde{w}(t_n') -  \tilde{w}_{\infty}|_H^2 
 +  \int_{t_n'}^{t_n'+2} E(\tau) d\tau. \label{eq230}
\end{equation}
By the assumption (A6) and Lemma \ref{w4} for any $\varepsilon > 0$ we can take a positive integer $N_1$ such that 
$$ \frac{1}{2}   |\tilde{w}(t_n') -  \tilde{w}_{\infty}|_H^2 < 
\frac{1}{2} \varepsilon \mbox{ for } n \geq N_1 \mbox{ and  }
E(t) < \frac{1}{4} \varepsilon \mbox{ for } t \geq N_1. $$
Then for $t \geq N_1 + 1$ there exists a positive integer $n \geq N_1$ such that $n + 1 \leq t \leq n+2$. In this case we see that 
$n \leq t_n' \leq n+1 \leq t \leq n + 2 \leq t_n'+2$. Hence, on account of 
(\ref{eq230}) we infer that 
$\frac{1}{2} |\tilde{w}(t) -  \tilde{w}_{\infty}|_H^2 < \varepsilon$ 
for $t \geq N_1 + 1$. This is the conclusion of this theorem. 
Thus we have proved 
Theorem \ref{large}. \hfill $\Box$ 





\bibliographystyle{alpha}

\end{document}